\documentclass{amsart}
\usepackage[T1]{fontenc}
\usepackage[english]{babel}
\usepackage{amsfonts,amsmath,amssymb,amsthm,enumitem}
\usepackage{graphicx,mathtools}
\usepackage{tikz}
\usetikzlibrary{backgrounds,calc}
\usepackage{floatrow}

\usepackage{cite}
\input xy
\xyoption{all}

\newcommand{\CC}{\mathbb{C}}

\newcommand{\HH}{\mathbb{H}}

\newcommand{\RR}{\mathbb{R}}

\newcommand{\ZZ}{\mathbb{Z}}

\newcommand{\tnu}{\tilde{N}}
\newcommand{\reg}{_{reg}}

\newcommand{\Lg}{\mbox{$\mathfrak g$}}

\newcommand\liegr{\sf}

\newcommand{\SU}[1]{\mbox{${\liegr SU}(#1)$}}

\newcommand{\U}[1]{\mbox{${\liegr U}(#1)$}}
\newcommand{\Sp}[1]{\mbox{${\liegr Sp}(#1)$}}
\newcommand{\SO}[1]{\mbox{${\liegr SO}(#1)$}}

\newcommand{\OG}[1]{\mbox{${\liegr O}(#1)$}}
\newcommand{\Spin}[1]{\mbox{${\liegr Spin}(#1)$}}
\newcommand{\G}{\mbox{${\liegr G}_2$}}
\newcommand{\F}{\mbox{${\liegr F}_4$}}
\newcommand{\E}[1]{\mbox{${\liegr E}_{#1}$}}
\newcommand{\T}[1]{\mbox{${\liegr T}^{#1}$}}

\newtheorem{cor}[equation]{Corollary}
\newtheorem{lem}[equation]{Lemma}
\newtheorem{prop}[equation]{Proposition}
\newtheorem*{thm*}{Theorem}
\newtheorem*{cor*}{Corollary}
\newtheorem*{lem*}{Lemma}
\newtheorem*{prop*}{Proposition}

\theoremstyle{definition}

\newtheorem*{defn*}{Definition}

\theoremstyle{remark}
\newtheorem{rem}[equation]{Remark}

\newtheorem*{rem*}{Remark}
\newtheorem*{ex*}{Example}
\newtheorem*{assumption*}{Assumption}

\newtheorem*{pf}{Proof}

\newcommand{\cref}[1]{Corollary~\ref{#1}}

\numberwithin{equation}{section}
\sffamily

\begin{document}

\title{Representations  of  compact  Lie  groups of\\
low cohomogeneity}
%{The Cohomogeneity of Compact Lie Group Representations }

\author[F.~J.~Gozzi]{Francisco J. Gozzi}

\address{Instituto de Matem\'atica e Estat\'\i stica, Universidade de
S\~ao Paulo, Rua do Mat\~ao, 1010, S\~ao Paulo, SP 05508-090, Brazil}
\email{fj.gozzi@gmail.com}

\thanks{The author has been supported by
  the FAPESP fellowship 2014/22568-1.}


\date{}%\today} 

\begin{abstract}
  We survey different tools to classify representations of compact Lie groups
  according to their cohomogeneity and apply these methods to the case
  of irreducible representations of cohomogeneity $6$, $7$ and $8$. 
\end{abstract}

\maketitle

\section*{Introduction}
%Compact connected Lie groups and their irreducible representations are algebraic in nature, 
%their classification via highest weight theory being well known.

%The study of linear representations of compact Lie groups is of fundamental importance 
%in regard to that of more general proper Lie group actions on manifolds and of the Alexandrov geometry of their orbit spaces.
%Hence our interest in low dimensional examples as a natural testing ground %in order
%to help us broaden our understanding of real irreducible representations 
%from the geometric viewpoint of their orbit spaces. % (cf. \cite{GL_OnOrbitSpaces}).

A fundamental invariant of a finite-dimensional
representation of a compact Lie group
is its \emph{cohomogeneity}, which by definition
is the minimal codimension $c$ of its orbits, and also agrees
with the topological dimension of the orbit space.

%It is apparent that r
Representations with low cohomogeneity should
display interesting geometric and topological properties. 
Indeed those with the minimum possible value $c=1$ 
are necessarily irreducible and coincide 
with the transitive isometric actions on unit spheres; 
they were probably first explicitly listed by Borel and Montgomery-Samelson
\cite{Borel,MS}, who in fact classified (effective) transitive smooth 
actions of compact connected Lie groups on spheres. 
The cases $c=2$ and $c=3$
were studied by Hsiang-Lawson \cite{HsiangLawson} 
(see also~\cite{str_ontheinvarianttheory}),
in connection with the construction of minimal submanifolds with large
groups of symmetries. They noticed that the maximal connected
groups with $c=2$ always act by the isotropy representation of a
symmetric space, and the same is true in the irreducible case
if $c=3$, up to three exceptions. More generally, representations
with $c=2$ are always polar~\cite{BCO}, and the irreducible ones with $c=3$
are always taut~\cite{GT}.

The complexity of representations grows with $c$ and it is
therefore natural to consider first irreducible representations. 
In \cite{GL_OnOrbitSpaces} irreducible representations with $c=4$ and
$c=5$ were classified 
in connection to certain problems regarding the geometry of orbit spaces.
Also, some interesting representations with $c=7$ associated to
quaternion-Kahler
symmetric spaces naturally appeared in~\cite{GG_QKreductions}. 
The purpose of this paper is to organize, 
extend and refine the techniques from~\cite{GL_OnOrbitSpaces} to
allow for a classification of representations of higher cohomogeneity.
In view of the inductive nature of our method, often requiring information
about lower cohomogeneity cases, 
we found useful to do a calculation of some principal isotropy groups.
These techniques are finally applied to classify irreducible
representations with $6\leq c\leq8$ and we obtain:

%Hence, the problem of computing a principal isotropy group can be addressed from 
%a Lie theoretic perspective and, as one may expect, this suffices to deal with the case of simple Lie groups altogether.
%In general we are required to consider arbitrarily big groups 
%and high degree representations which, in all, have small quotients.
%Our focus is then set on lower cohomogeneity estimates and on working with tensor products
%since the difficulty lays in non-simple Lie groups.

\begin{thm*}
A compact connected Lie group irreducible representation of cohomogeneity 6, 7 or 8 is either polar or, 
else, listed in Tables \ref{table cohom 6}, \ref{table cohom 7} or \ref{table cohom 8}.
\end{thm*}

%It should be possible to extend the classification up to higher cohomogeneity since
%the general procedure has an inductive flavor. 
%In doing so, one requires more information about the previous lower cohomogeneity cases, 
%hence our interest in determining the principal isotropy groups whenever possible.

This paper is organized as follows. 
The first section establishes terminology and preliminary considerations, 
where useful principal isotropy computations and cohomogeneity estimates are 
given for sums and tensors of representations.
The second section represents the core of this work where the available tools are discussed, 
distinguishing among simple and non-simple Lie groups.
%Here we also summarize and rectify a result of \cite{Hsiang_Diff_II} concerning the principal isotropy of 
%compact simple Lie group representations.
Last, the third section sketches the classification for arbitrary cohomogeneity
and illustrates the procedure by carrying it through up to cohomogeneity $8$. 
The main tables are displayed at the end. 

The author is indebted to Claudio Gorodski for suggesting this project and for many fruitful discussions, his was the suggestion for the proof of
Remark \ref{rem Claudio}.

\section{Preliminaries}\label{prelim}

Let \mbox{$\rho: G \to \OG V$} be an orthogonal representation of a
compact Lie group $G$ on a finite-dimensional Euclidean vector space~$V$. 
The stabilizer of a point $p\in V$ is denoted by $G_p$ and referred to as its \textit{isotropy group}, and the orbit through $p$ is denoted by $Gp$. 
Isotropy subgroups along a $G$-orbit lie in the same conjugacy class.
The conjugacy classes of isotropy groups, called \textit{isotropy types},
give a natural stratification of the orbit space. 
Moreover, the inclusion of subgroups induces a lattice structure on the
set of isotropy types. The Principal Orbit Type Theorem asserts
that there is a unique minimal istropy type, called the
principal isotropy type; the corresponding orbits are called
\emph{principal orbits}. The union of principal orbits is a
open dense subset of $V$, and its projection to the orbit
space is a connected, open dense subset called the \emph{principal stratum}.
We will denote a fixed principal isotorpy group
by~$H$.  

The cohomogeneity of a representation equals the codimension of a
principal orbit $G/H$ in $V$, 
and thus it can  be computed from $(H)$ from the following equation:
\begin{equation}
c=c(\rho) = \dim (V) - \dim(G) + \dim(H). 
\end{equation}
%where we write $\rho=:\!(G,V)$ when the action is clear from the context.

%\subsection{Naive Algorithm}\label{subsection naive algorithm}
There is an algorithm to compute a principal isotropy group of a
representation $\rho$ based on the fact that the set
of isotropy groups of any the slice
representation of $\rho$ contain representatives of all
isotropy types of $\rho$~\cite{hsiang-smith-theory}.
Namely, take any point $p$ in $V$ which is not fixed by $G$.
Then there is a $G_p$-invariant
decomposition 
\[ V= T_p (Gp)\oplus N_p(Gp) \]
where $T_p(Gp)$ and $N_p(Gp)$ denote the tangent and normal spaces
to $Gp$ at $p$, respectively. 
The \emph{slice representation} of $\rho$ at $p$ is 
the action of $G_p$ on $N_p(Gp)$. A principal isotropy group for
$(G_p,N_p(Gp))$ is a principal isotropy group for $(G,V)$ as well, 
and we have reduced our calculation to the case of a representation of lower
dimensional group and lower dimensional representation space. By dimensional
reasons, the algorithm evetually stops when the slice representation is
trivial.

In general, the calculation of full isotropy groups can be a delicate matter
and in most cases we will only need their isotropy groups which are much easier
to obtain. It is useful to note that, by considering only the identity
component of the isotropy group at each step of the algorithm, we 
end up with the identity component of the principal isotropy group.
%Notice that for a complex representation it is natural to perform this computations on the level of Lie algebras.
%If the representation is given in terms of a maximal weight of $\Lg$, 
%we may start with a vector associated to the maximal weight and easily determine its isotropy which, to begin, contains all positive roots. 
%The tangent space to the orbit follows easily as well. We remark that for a real type representation one can proceed analogously, 
%as long as we consider only those vectors which are invariant under the real structure automorphism.

\subsection{Reducible representations}\label{subsection reducible}
Although we are mostly interested in the case of irreducible representations,
reducible ones are bound to appear for instance 
as slice representations at singular points. It is plain that
for a $G$-invariant decomposition $V=V_1\oplus V_2$, a point
$p=(p_1,p_2)$ with $p_i\in V_i$ has isotropy group $G_p=G_{p_1}\cap G_{p_2}$;
and $G_{p_i}$ 
being a principal isotropy group of $(G,V_i)$ for $i=1$, $2$ is a necessary
but not sufficient condition for $G_p$ to be a principal isotropy group of
$(G,V)$. In general, it follows from the discussion above about slice
representations that a principal isotropy group of $(G,V)$ coincides
with a principal isotropy group of $(H_1,V_2)$, where $H_1$ is a
principal isotropy group of $(G,V_1)$.

For example, given two copies of a non-trivial representation $(G,V)$,
note that a principal isotropy group $H$ acts reducibly on $V$,
where $V^{H}$ is a proper subspace of dimension equal to $c(G,V)$, and that 
$c(G,2V) = c(G,V) + c(H,V)$ and $c(H,V)	= c(G,V) +  c(H, T_{[H]} G/H)$, 
 so 
 \begin{equation}\label{eq cohom of m copies}
c(G,2V)= 2c(G,V) + c(H, T_{[H]} G/H)\geq 2c(G,V)+ 1.
 \end{equation}
% For $m$ copies, we have a rough lower bound $c(G, mV) \geq m \cdot c(G,V) + m -1$ 
% though this may be substantially improved .

\subsection{Real tensor products}\label{subsection real tensors}
%\noindent\textbf{Real Tensors}

Let $\rho:G\to \OG V$ be a representation of the form
$\rho=\rho_1\otimes\rho_2$, where $G=G_1\times G_2$, $V=V_1\otimes_{\mathbb R}V_2$
and $\rho_i:V_i\to \OG{V_i}$ for $i=1$, $2$. We consider
the slice representation at a pure tensor $p=p_1\otimes p_2$,
where $H_i=G_{p_i}$ is a principal isotropy group for $i=1$, $2$. 

Note first that, for the identity components, $G_p^0=H_1^0\times H_2^0$. 
Furthermore, since the Lie algebra of $G$ acts on $V_1\otimes V_2$ by
derivations, the tangent space
\[T_p(Gp) = 
\left( T_{p_1} (G_1p_1) \otimes p_2  \right)\oplus \left( p_1 \otimes T_{p_2} (G_2p_2) \right).\] 
The slice representation has $\mathbb R p$ as a trivial
summand, so it is natural to consider the \textit{reduced normal space}
$\tnu_p(Gp)= N_p (G\cdot v)\ominus\mathbb R p $, 
and then
\begin{equation}\label{tilde-N}
 \tnu_p(Gp) = 
 (\tnu_{p_1}(G_1p_1) \otimes p_2 ) \oplus (p_1 \otimes \tnu_{p_2}(G_2p_2)) 
 \oplus (p_1^{\perp} \otimes p_2^{\perp}).
 \end{equation}
 The group $H_1^0\times H_2^0$ acts trivially on the first two summands of
 the right hand-side of~(\ref{tilde-N}) 
leading to 
\begin{equation}\label{eq real tensor}
c(G,V) = c(G_1, V_1) + c(G_2,V_2) - 1 + c(H_1^0\times H_2^0, p_1^{\perp} \otimes  p_2^{\perp}) .
\end{equation}

It is useful to discuss further
the last term on the right hand-side of~(\ref{eq real tensor})
In some particular cases.
%A further study of the last term proves useful in particular cases. 
Write $T_{p_i}(G_ip_i)=U_i$, $\tnu_{p_i}(Gp_i)=\tnu_i$ and $c(G_i,V_i)=c_i$.
%We have the following direct sum decomposition,
%$$v_1^{\perp} \otimes  v_2^{\perp}=[\tnu_1\otimes \tnu_2] \oplus [\tnu_1\otimes U_2] \oplus [U_1\otimes U_2] \oplus [U_1\otimes \tnu_2] .$$
If $c_1=1$, then 
$v_1^{\perp} \otimes  v_2^{\perp}$ equals $ (U_1\otimes U_2) \oplus (U_1\otimes \tnu_2)$.
The action of $H_1^0\times H_2^0$ on the second summand equals to $c_2-1$ copies of $(H_1,U_1)$, whose principal isotropy group is therefore
of the form $K_1\times H_2^0$ for a subgroup $K_1$ of $H_1^0$. 
We finally obtain that 
\begin{equation}\label{eq real tensor cohomo 1}
c(G,V) = c_2 + c(H_1^0,(c_2-1)U_1) + c(K_1\times H_2^0, U_1\otimes U_2).
\end{equation}

\iffalse
Finally, let us remark that, in the general plan of enumerating representations by their cohomogeneity,
those with lower values of $c$ continue to appear as factors of tensor products of higher cohomogeneity.
Hence, in the light of equation \eqref{eq real tensor} we emphasize the usefulness of determining 
principal isotropy types or at least their identity connected component.
\fi

\medskip

\section{Main}

\subsection{Simple Lie Groups}\label{section simple lie groups}

%by the work of A. Kollross in \cite{Kollross_HyperpolarCohom1} and \cite{Kollross_LowCohomAndPolarActions}, 
A representation of a compact simple Lie group $G$ with a certain bound
on the cohomogeneity either satisfies the condition of 
\textit{low degree} defined by A. Kollross
in~\cite{Kollross_HyperpolarCohom1} which leads to a
concrete list of candidates,
or $G$ has bounded dimension.
In the latter case, we can resort 
to a case by case analysis of the finitely many possible groups 
making use of basic representation theory. For a fixed maximal
torus of $G$ and an ordering of the roots.  
Then we have the simple roots by $\alpha_1,\ldots,\alpha_r$ and 
the fundamental weights $\lambda_1,\ldots,\lambda_r$, and 
the degree %or complex dimension 
of a complex irreducible representation is a
strongly monotonous function of the $\lambda_i$, 
as observed by A. L. Onishchik in \cite{Onishchik}.
%\textcolor{red}{So is claimed by Kollross in ``A classification of hyperpolar and cohomogeneity one actions''}.
By making use of the basic inequality
$$ c(\rho)\geq\mathrm{deg}(\rho) - \dim (G),$$ %deg(\rho(c_1,\ldots,c_r)) - dim (G) $$
we may therefore bound the degree for each possible group and obtain
a finite list of representations.

\smallskip
The problem of determining the principal isotropy type and, hence, 
the exact cohomogeneity of a given $n$-dimensional representation is
completely solved in the case of compact simple Lie groups
in~\cite{Hsiang_Diff_II}.
%In fact, it should be noticed that the former includes the case of discrete (non-trivial) principal isotropy groups as well.
%For convenience we state a rectified version of their classification in the following.
It is interesting to reformulate their result as follows.

\begin{prop}%[W. C. Hsiang \& W. Y. Hsiang]
\label{prop simple pif Hsiang}
An irreducible representation of a compact connected simple Lie group
with non-trivial principal isotropy type is either polar or
a simple factor of the isotropy representation associated to a symmetric space
of Hermitian or quaternion-K\"ahler type.
\end{prop}

We recall that a representation is called \emph{polar}
if there exists a subspace, called a \emph{section},
meeting all orbits and always orthogonally. It follows from
the classification of Dadok~\cite{Dadok} that a polar representation
is orbit-equivalent to (i.e.~has the
same orbits as) the isotropy representation of a symmetric space,
and indeed the maximal groups in each orbit-equivalence class
act by isotropy representations of symmetric spaces. 
It turns out that the only irreducible polar representations with trivial
principal isotropy groups are the standard actions of $\SO2$ on $\mathbb R^2$
and of $\SU2$ on $\mathbb C^2$~\cite{GL3}.  

A representation $(G,V)$ with non-trivial principal isotropy group $H$
admits the so called \emph{core reduction}, or \emph{Luna-Richardson-Straume reduction} in the linear case, namely,
the quotient group $N_G(H)/H$ acts effectively
on the fixed point set $V^H$ with orbit space
isometric to the orbit space $V/G$. More generally, a \emph{reduction}
of $(G,V)$ is a representation $(K,W)$ with $\dim K<\dim G$ and
orbit space isometric to $V/G$~\cite{GL_OnOrbitSpaces}. 
In the special case $N_G(H)/H$ is finite,
the representation $(G,V)$ is called \textit{asystatic}, 
and it is automatically polar with $V^H$ as a section 
(\cite{Alekseevsky_asystatic}, see also~\cite{GroveZiller_polar}). 
Remarkably, every polar representations admits an orbit-equivalent
finite extension to an asystatic representation of a
disconnected Lie group, as observed in
\cite{str_ontheinvarianttheory} (see also~\cite{KollrossGorodski}).
The non-polar representations satisfying
the conditions in Proposition~\ref{prop simple pif Hsiang}
are given by the restriction to the non-$\U 1$ factor
(resp.~non-$\Sp 1$ factor)
of the isotropy representation of Hermitian (resp.~quaternion-K\"ahler
symmetric space).
They admit reductions to tori (resp.~$\Sp 1^3$-subgroups),
and were investigated in~\cite{GL2} and~\cite{GG_QKreductions}.

\iffalse

The original work in \cite{Hsiang_Diff_II} included the case of a half-spin 
representation of $\Spin {14}$ of complex type and degree $64$, 
which was incorrectly claimed to have non-trivial discrete principal isotropy. %, see Subsection \ref{subsec reductions}.
In fact, one may use the method of taking repeated slices to verify that its principal isotropy is indeed trivial 
-though one must do so carefully in order to determine the full principal isotropy and not just its connected component.
Alternatively, one can use the methods of \cite{GL_OnOrbitSpaces,GG_QKreductions} to check that its orbit space has trivial boundary, 
which also implies the claim.

\begin{rem} \label{rem simple list}
%For reference, 
A list of compact simple Lie group representations with $2\leq c \leq 8$ can be found in \cite{GL_OnOrbitSpaces}. 
This is easily enlarged to contemplate $c=9$ 
by extending parameters for various families of examples listed therein
plus one more case given by a degree $32$ representation of $\Spin {11}$ of quaternionic type. 
%Such a list is a necessary input in order to classify irreducible representations up to cohomogeneity $8$ 
%in Section \ref{section classification}. 
\end{rem}

\fi

We shall refer to the Tables in Subsection~12.8 in \cite{GL_OnOrbitSpaces}	 listing polar and non-polar
irreducible representations of compact connected simple Lie groups of
cohomogeneity up to $8$. In addition, we extend those lists to 
include cohomogeneity $9$ in Table~\ref{table cohom 9}. 

\begin{table}[h]\label{table cohom 9}
\[ \begin{array}{|ccc|}
\hline
G & V &  \textrm{Polar?} \\ 
\hline
\SO{10} & S_0^2\RR^{10} & \mbox{Yes}\\
\Sp{10} & \Lambda^2\CC^{20}\ominus\CC & \mbox{Yes} \\
\SU{10} & \mbox{Adjoint} & \mbox{Yes} \\
\SO{18} & \mbox{Adjoint} & \mbox{Yes} \\
\SO{19} & \mbox{Adjoint} & \mbox{Yes} \\
\Sp{10} & \mbox{Adjoint} & \mbox{Yes} \\
\SU{19} & \mbox{Adjoint} & \mbox{Yes} \\
\SU2 & \CC^6 & \mbox{No} \\
\SU8 & S^2\CC^8 & \mbox{No} \\
\Spin{11} & \HH^{16} & \mbox{No}\\
\hline
\end{array}\]
\caption{Irreducible representations of cohomogeneity~$9$ of 
simple Lie groups.}
\end{table}

\begin{rem} \label{rem Claudio}
  At this juncture, it is worth commenting on the half-spin representation
  of $\Spin{14}$ on $\mathbb C^{64}$ that appears on~\cite[Table~A]{Hsiang_Diff_II}.
  It can be shown to have trivial principal isotropy group by an argument
  similar to  Example~1 in that paper. Indeed, if $p$ is a highest weight vector,
  it is easy to see that $G_p^0=\SU7$. Moreover, from the theory
  of parabolic subgroups, it is known that the normalizer of the line
  $\mathbb R p$ is $\U7$. 
  Finally, let $T$ and $\tilde T$ be maximal tori in $\SU7$, $\U7$,
  respectively. 
  It is easy to see that no non-trivial elements of $\tilde T/T$ fix~$p$ and, thus, 
  $G_p=\SU7$ is the full isotropy subgroup. Furthermore, its slice representation is given by $\RR \oplus \CC^7\oplus \Lambda^3\CC^7$ 
  which has trivial principal isotropy.
  \end{rem}

\subsection{Non-simple Lie groups}

For the case of real representations of
non-simple compact Lie groups, it is easier to
start with complex representations 
since they have a simple description
in terms of tensor products of complex irreducible
representations of the factors of $G$. 

A complex representation $\pi$ is called of \emph{real type}
if it comes from a representation on a real vector space by extension
of scalars, and it is called of \emph{quaternionic type} if it comes
from a representation on a  quaternionic (right-)vector space by
restriction of scalars. If $\pi$ is neither of real type nor of
quaternionic type, we say that $\pi$ is of \emph{complex type}.
If $\pi$ is irreducible, then it is exactly of one of those types. 

Now it is known that the finite-dimensional 
real irreducible representations $\rho$ of $G$ 
fall into one of the following disjoint classes:
\begin{enumerate}
\item[(a)] the complexification $\rho^c$ is irreducible and $\rho^c=\pi$ is 
a complex representation of real type;
\item[(b)] the complexification $\rho^c$ is reducible and 
$\rho^c=\pi\oplus\pi$ where $\pi$ is a complex irreducible representation of 
quaternionic type;
\item[(c)] the complexification $\rho^c$ is reducible and 
$\rho^c=\pi\oplus\pi^*$ where $\pi$ is a complex irreducible representation of 
complex type and $\pi^*$ is not equivalent to $\pi$ 
(where $\pi^*$ denotes the dual representation of $\pi$).
\end{enumerate}
The relation between $\rho$ and $\pi$ is that $\rho$ is a real form 
of $\pi$ in the first case ($\rho^c=\pi$ and $\rho=[\pi]_{\mathbb R}$). In the other two cases
$\rho$ is $\pi$ viewed as a real
representation, its ``realification''. 
We shall call $\rho$ of \emph{real}, \emph{quaternionic} or \emph{complex type} according to whether the associated 
$\pi$ is of real, quaternionic or complex type.
Note also that $\pi$ is self-dual precisely in the first two cases.

\iffalse

Our focus is set on real irreducible representations which give rise, by complexification, 
to unique complex representations. Reciprocally, 
an irreducible complex representation may still be irreducible over the real number field or, else, 
admit a totally real irreducible invariant subspace $W$ such that $V=W\otimes_{\RR}\CC$.
In the former case the representation is of complex or quaternionic type 
while in the latter the (complex) representation $V$ is of real type and $[V]_{\RR}$ denotes a choice of a real form, 
unique up to isomorphism.

\fi

A real (resp.~quaternionic) structure for $\pi$ can be equivalently
given by an equivariant conjugate-linear endomorphism
$\epsilon$ such that $\epsilon^2=1$ (resp.~$\epsilon^2=-1$).
It follows that if $\epsilon_i$ is such a structure for $\pi_i$ then
$\epsilon_1\otimes\epsilon_2$ is such a structure for
$\pi_1\otimes_{\CC}\pi_2$, and $\pi_1\otimes_{\mathbb C} \pi_2$ is
of complex type whenever one of the factor is.
Table \ref{table type of tensor repr} summarizes our discussion.

\begin{table}[h!]
\caption{
%Double entry chart showing 
The type of $\rho_1 \otimes_{\CC}\rho_2$ in terms of $\rho_1$ and $\rho_2$.}
\label{table type of tensor repr}
\begin{minipage}[r]{5cm}
$$\begin{array}{c|ccc}
	&	r	&	q	&	c	\\
\hline
 r	&	r	&	.	&	.	\\
 q	&	q	&	r	&	.	\\
 c	&	c	&	c	&	c 	\\
\end{array}$$
\end{minipage}
\end{table}

\begin{rem} \label{rem tips for tensors by types}
We collect some useful remarks
about real forms and realifications of tensor products of complex representations
that will allow us to use ideas from Section~\ref{subsection real tensors}.
\iffalse

One may always consider the associated complex representation first, determine its type and, accordingly, find a real form or 
consider the underlying real structure.
However, we prefer to work with real tensors whenever possible 
%since we have a deeper insight for computing slices, principal isotropy groups and, hence, cohomogeneity, 
and have the results of Section \ref{subsection real tensors} at our disposal.
We collect some useful observations.
\fi
\begin{enumerate}[label=(\roman*)]
\item If $\pi_1$ and $\pi_2$ are of are of real type,
  a real form of $\pi=\pi\otimes_{\CC}\pi_2$ can be obtained as
  $[\pi]_{\RR}= [\pi_1]_{\RR} \otimes_{\RR} [\pi_2]_{\RR} $.
\item If $\pi_1$  is of real type and $\pi_2$ is of complex or
  quaternionic type, the realification of $\pi=\pi_1\otimes_{\CC}\pi_2$
can be described as $\pi^r = [\pi_1]_{\RR} \otimes_{\RR} \pi_2^r$. 
\item If $\pi_1$ and $\pi_2$ are of quaternionic type,
  a real form of $\pi=\pi_1\otimes_{\CC}\pi_2$ can be obtained as
  $[\pi]_{\RR}= \pi_1 \otimes_{\HH} \pi_2$.
  We refer to~\cite[Subsection~12.3]{GL_OnOrbitSpaces} for
an explanation about quaternionic tensor products. 
\item When $\pi_1$ is of complex type and $\pi_2$ is not of real type,
  we must deal directly with the realification of
  $\pi_1\otimes_{\CC}\pi_2$.
\end{enumerate}
\end{rem}

%\subsection{Tensors of cohomogeneity one representations} \label{section Tensors of cohomogeneity one representations}

A representation %$(G,V)$ of a non-simple Lie group 
which is given as a real, complex or quaternionic tensor product
is naturally a restriction of 
$(\SO m \times \SO n, \RR^m \otimes_{\RR} \RR^n)$, $(\U m \times \U n, \CC^m \otimes_{\CC} \CC^n)$ or
$(\Sp m \times \Sp n, \HH^m \otimes_{\HH} \HH^n)$, respectively. 
Hence, we obtain a lower bound for the cohomogeneity, namely,
the minimum between $m$ and $n$. 
%Notice that these representations are polar and the 
For simplicity, below we only list the
identity components of the corresponding principal isotropy groups
(see for instance \cite{Hsiang_Diff_II}).

\begin{lem}%[W. C. Hsiang \& W. Y. Hsiang]
For a representation $(G,V)$, let $H$ denote a principal isotropy group.
Assume $m\geq n$. Then:
\begin{enumerate}
 \item For $(\SO m \times \SO n, \RR^m \otimes_{\RR} \RR^n)$,  %$\rho=\SO m \otimes_{\RR} \SO n$, 
 we have 	$H^0=\SO{m-n}$.
 \item For $(\U m \times \U n, \CC^m \otimes_{\CC} \CC^n)$,
 %$\rho=\U {m}\otimes_{\CC} \U {n}$, 
 we have 	$H^0=\U{m-n}\times \U 1^n$.
 \item For $(\Sp m \times \Sp n, \HH^m \otimes_{\HH} \HH^n)$, we have
 %$\rho=\Sp m\otimes_{\HH} \Sp n$, 
  $H^0=\Sp {m-n}\times \Sp 1^n$.
 \end{enumerate}
%The respective subgroups being block-embedded into the first factor.
 \end{lem}

\iffalse
Moreover, since our exposition is organized according to the cohomogeneity of the factors of a tensor representation,
the case where both factors are of cohomogeneity one turns out to be a natural first step.
Motivated by this and the previous Lemma, we consider real or complex tensor representations among the standard
$\SO n$, $\U n$ and $\Sp n$ vector representations. We have the following.
\fi

With a view toward complex tensor products of mixed type,
we state the following result whose proof is straightforward based on the algorithm described in Section~\ref{prelim}.
From the information about principal isotropy groups we can of course
deduce the cohomogeneity of the involved representations. 

\begin{lem}\label{lem Tensors of cohom 1 representations}
For a representation $(G,V)$, let $H$ denote a principal isotropy group.
Then:

\begin{enumerate}
%1
\item For $(\SO m \times \U {n}, \RR^m\otimes_{\RR} \CC^n)$ we have
 \[H^0= \begin{cases}
\SO{m-2n}& m > 2n + 1, \\
\U{n-m}  & n > m,
\end{cases}\]
otherwise it is trivial.
%2
\item For $(\SO m \times \Sp {n}, \RR^m\otimes_{\RR} \HH^n)$, we have 
\[ H^0= \begin{cases}
\SO{m-4n}&\mbox{if $m > 4n + 1\geq5$;} \\
\Sp{n-m}  &\mbox{if $n > m\geq3$;}\\
\end{cases}\]
otherwise it is trivial.
%3
\item For $(\SO m \times \Sp 1 \Sp {n}, \RR^m\otimes_{\RR} \RR^{4n})$, 
\[H^0= \begin{cases}
\SO{m-4n}	& m > 4n + 1, \\
\Sp{n-m} 	& n > m\geq 3,
\end{cases}\]
otherwise it is trivial, as long as $m\geq 3$.
For the case $m=2$ and $n\geq 2$ we have that $H^0= \{1\}\times \U 1 \Sp {n-2}$. 
%4
\item For $(\SO m \times \U 1\Sp n, \RR^m\otimes_{\RR} \CC^{2n})$ we have
\[H^0= \begin{cases}
\SO{m-4n}& m > 4n , \\
\Sp{n-m}  & n \geq m,
\end{cases}\]
otherwise it is trivial. 
%5
\item For $(\U m \times \Sp 1\Sp n, \CC^m\otimes_{\RR} \RR^{4n})$ we have
\[H^0= \begin{cases}
\U{m-4n}& m \geq 4n , \\
\Sp{n-2m}  & n \geq 2m,
\end{cases}\]
otherwise it is trivial, assuming $n\geq 2$. 
%6
\item For $(\Sp 1\Sp m \times \Sp 1\Sp n, \RR^{4m}\otimes_{\RR} \RR^{4n})$ we have 
%For $\SO m \otimes_{\RR} \U {n}$ the Id-component of a principal isotropy is given by 
\[H^0 = \Sp{m-4n}\]
if $m > 4n + 1$, otherwise it is trivial. 
%7
\item For $(\U {m}\times \Sp n, \CC^m\otimes_{\CC} \HH^n)$, we have
\[ H^0= \begin{cases}
\U {m-2n}	&\mbox{if $m\geq n+1\geq2$;} \\
\Sp {n-m}	&\mbox{if $n\geq m\geq3$;}\\
\U1\times\U {m-2}	&\mbox{if $m\geq2$ and $n=1$;} \\
\U1\times\Sp{n-m} &\mbox{if $n\geq m$ and $2\geq m\geq1$;}\\
\end{cases}\] 
otherwise it is trivial. 
%8
\item For $(\Sp {m}\times \U 1\Sp n, \HH^m\otimes_{\CC} \CC^{2n})$, we can assume $m\geq n\geq 2$, then 
$ H^0= \Sp {m-2n}$ if $m > 2n$ or, else, it is trivial. 
\end{enumerate}
\end{lem}

%Analogous results can be obtained for tensors between classical families of cohomogeneity one representations.

For each family of representations listed in
Lemma~\ref{lem Tensors of cohom 1 representations} with fixed $m$ (resp.~$n$),
one notes that the cohomogeneity is
asymptotically constant on $n$ (resp.~$m$).
This can be deduced in general from the following useful result.

\begin{lem}[Monotonicity Lemma\cite{GL_OnOrbitSpaces}]\label{monoton}
  Let $\rho(n)$ be the $\mathbb F$-tensor product of a fixed real
  representation of $\mathbb F$-type with $\SO n$, $\U n$ or $\Sp n$
  according to whether $\mathbb F=\RR$, $\CC$ or $\HH$.
  Then the cohomogeneity of $\rho(n)$ is a non-decreasing
  function of $n$.
\end{lem}

%\begin{lem}
%  The cohomogeneity of the representation~$\rho(n)$
%  given in the statement of Lemma~\ref{monoton} is
%  asymptotically constant on $n$. 
%\end{lem}

\begin{cor}
The cohomogeneity of any given representation tensored with one of the standard
vector representations: $(\SO n,\RR^n)$, $(\U n,\CC^n)$ or $(\Sp n,\HH^n)$, is asymptotically constant on $n$. 
\end{cor}

\begin{pf}
  In case $\mathbb F=\RR$, assume the real dimension of a given
  representation is $k$. For $n>k$, the principal isotropy group
  of $\rho(n)$ contains a normal subgroup of type $\SO{n-k}$,
  which yields an upper bound on the cohomogeneity. Lemma~\ref{monoton}
  then implies that the cohomogeneity stabilizes, as wished.
The other cases are analogous.  
\qed
\end{pf}

\medskip

 \section{Classification} \label{section classification}

Let us sketch a general plan for the classification of 
irreducible representations of compact connected Lie group 
with bounded cohomogeneity. We shall then carry this procedure up to cohomogeneity $c \leq 8$.

\smallskip

First, as discussed in Section \ref{section simple lie groups},
we may determine the irreducible representations of simple Lie groups up to a given cohomogeneity.
This task is also an auxiliary one for the cases to come and, hence, we shall see that it suffices to list them
up to cohomogeneity $c+3$. % when looking for a general classification up to cohomogeneity $c$.  
%In fact, one may do better in specific computations and, in our classification up to cohomogeneity $8$ we shall only need to list them up to $c \leq 9$. 

On a second step, we consider those representations that can be expressed as a real tensor product, see \mbox{Remark \ref{rem tips for tensors by types}}. 
Though a great deal of cases is to be considered here we have the useful formula \eqref{eq real tensor} at our disposal. 
Observe that both the two factors involved in the tensor have strictly smaller cohomogeneity and are, thus, inductively known. 
%The computation of the exact cohomogeneity is aided by any knowledge of the previous principal isotropy types of the factors.

Third, we have the case of the quaternionic tensor product of two irreducible representations of quaternionic type.
%Observe that each of quaternionic type factor either corresponds to a simple Lie group or, else, 
%is a product of a real type representation times a quaternionic one.
Notice that we may assume the group to be a product of two simple factors for, otherwise, 
we would be able to express the outcome representation as a suitable real tensor product, already considered. 
Moreover, a standard application of the Monotonicity Lemma \ref{monoton} shows that each (simple) factor has cohomogeneity at most $c+3$.
In the following classification up to cohomogeneity $8$ we consider more cases in this step and, as a result, 
we need only list simple Lie group representations up to cohomogeneity $9$.

Finally, the remaining cases are real representations underlying complex tensor products, 
again after Remark \ref{rem tips for tensors by types}. 
These representations split completely as a complex tensor product of its simple factors and, 
thus, give a subrepresentation of $\U {n_1} \otimes \SU {n_2} \otimes \cdots \SU {n_k}$. 
Monotonicity %(Lemma \ref{monoton}) 
can be used to give an upper bound on the number ``$k$'' of factors
by estimating the cohomogeneity of the initial case %a tensor product of standard vector representations of the form
%$\U 2 \otimes_{\CC} \SU{3} \otimes_{\CC} \cdots \otimes_{\CC} \SU{3}$
$(n_1,n_2,\ldots,n_k)=(2,3,\ldots,3)$. Here, once more, we avoid the case of two quaternionic $\SU 2$ factors which leads to a reducible representation. 
%we have at most one quaternionic factor. %, implying $3\leq m_2\leq m_k$.
Furthermore, each simple factor has cohomogeneity at most $c+1$.
%$\SU{m}\otimes_{\CC} \SU n\otimes_{\CC} \U p$ with $m\geq 2$ and $n,p\geq 3$, 
%%%%One may ``permute'' the circle $\U 1$ among the factors, and thus gain 
%Furthermore, we have monotonicity with respect to each of the variables: $m,n,p$, 
%for the cohomogeneity of the corresponding representation. 
%We then compute $c(\SU{m}\otimes_{\CC} \SU n\otimes_{\CC} \U p)\geq 16$, 
%from a natural bound for the initial case \mbox{$(m,n,p)=(2,3,3)$}. 

%\bigskip 
%
%\textbf{Classification up to cohomogeneity $8$}

We emphasize that this is an inductive procedure,
though at each step one is required to compute or estimate the cohomogeneity 
of finite many (possibly infinite) families of representations.

\medskip

We illustrate the procedure in carrying out the classification of irreducible
representations of compact connected Lie groups of cohomogeneity~$6$, $7$ or $8$. 
 
 \subsection{Real tensor products}
 Let us consider $\rho=\rho_1\otimes_{\RR} \rho_2$ as a
 restriction of the representation $\SO m\otimes_{\RR} \SO n$,
 where $\rho_1$ and $\rho_2$ are real irreducible representations,
 and assume $6\leq c(\rho)\leq 8$.
\subsubsection{}
Assume first $m=2$ so that $\rho_1=(\SO 2,\RR^2)$. 
Then the second factor $\rho_2$ must be of real type for otherwise
$\rho$ would be reducible.
Formula~\eqref{eq real tensor cohomo 1} gives
\begin{equation*}%\label{1}
 c(\rho)= 2c(\rho_2) -1 + c(H_2,U_2)\geq 2c(\rho_2). 
\end{equation*}
Now $c(\rho)\leq 8$ implies that $c(\rho_2)\leq 4$. 
Let us run through the possible cases.

First, we cannot have $c(\rho_2)=4$ for otherwise
$c(H_2,U_2)\geq3$ which is impossible.
In fact, this is clear if $\rho_2$ is a polar representation
since in this case its principal orbits are irreducible
isoparametric submanifolds admitting at least three
curvature distributions~\cite{BCO} which provides an $H_2$-invariant
decomposition of~$U_2$ with at least three components.
If $\rho_2$ is non-polar, then by~\cite[Table~1]{GL_OnOrbitSpaces}
it is $(\SO3,\RR^7)$ or $(\SO3\times\G,\RR^3\otimes_{\RR}\RR^7)$;
in both cases $H_2$ is finite and $c(H_2,U_2)\geq3$ follows. 

Second, the case $c(\rho_2)=3$ is discarded as well, 
by improving the previous estimate to $c(H_2,U_2)\geq 5$.
%Consider the possibilities for $\rho_2$.
In fact, 
%if $\rho_2$ is a polar representation, then using 
%again that a principal orbit $G_2p_2$ is an irreducible
%isoparametric submanifold, $c(H_2,U_2)<5$ can happen only if 
%$G_2p_2$ has at most $4$ principal curvatures. From general theory,
%this is the case only if $G_2p_2$ is of $\sf A_2$-type
%or $\sf B_2$-type; but in such cases, the codimension of $G_2p_2$
%is $2$ contradicting $c(\rho_2)=3$. On the other hand,
if $\rho_2$ is not polar, by classification we know 
that $\rho_2=(\Sp 1\times\Sp n, \HH^2\otimes_{\HH}\HH^n)$~\cite{HsiangLawson},
and it follows from results in~\cite{GT02} that $c(H_2,U_2)\geq 9$.
If not, polar candidates can be discarded by estimating the value of $c(H_2,U_2)$ directly, 
since their principal isotropy algebras are known. 
In the case of families of s-representations associated to the orthogonal and quaternionic grassmannians
it is enough to estimate the cohomogeneity of the initial cases and apply the Monotonicity Lemma.

Third, if $c(\rho_2)=2$ then $\rho_2$ is polar and $c(H_2,U_2)\leq 5$.
A principal orbit is isoparametric and we use their classification
to list the possibilities in Table \ref{table cohom 2 real type repr}; 
we obtain examples corresponding to the first four cases therein. 

\begin{table}[h]
$$ \begin{array}{|cc|ccc|}
\hline
G_2	&	V_2		&	H_2	& U_2		&	c(H_2,U_2)	\\
\hline
\hline
\SO 3			&\RR^5				&\{1\}		& 3\RR					&	3 	\\
\Sp 3			&[\Lambda^2\CC^6\ominus\CC]_{\RR}	& \Sp 1^3	& 3\HH					&	3	\\
\SU 3			&\mathfrak{su}(3)\cong \RR^8	& \T2		& 3\CC					&	4	\\
\F			&\RR^{26}			& \Spin 8 	&\RR^8_0\oplus \RR^8_+ \oplus \RR^8_-	&	4	\\
\Sp 2			&\mathfrak{sp}(2)\cong \RR^{10}	& \T 2		&4\CC					&	6	\\
\Sp 2\times \Sp n	&\HH^2\otimes_{\HH}\HH^n	&\Sp1^2 \Sp{n-2}& 2[\HH^{n-2}\oplus \RR^3\oplus \HH]	&> 6	\\
\SO4 & \RR^8 &  \ZZ_2^2 & 6\RR & 6 \\
\G & \Lg_2\cong\RR^{14} & T^2 & 6\CC & 10\\
\hline
\end{array}$$
 \caption{Isotropy representations $(G_2,V_2)$ of real type irreducible isoparametric submanifolds
of codimension~$2$. }\label{table cohom 2 real type repr}
\end{table}

Finally, if $c(\rho_2)=1$, then we recall $\rho_2$ if of real type
and refer to the classification to see that the outcome $\rho$ is either polar 
with cohomogeneity $2$
($G_2=\SO n$, $\G$ or $\Spin7$) or has cohomogeneity~$3$ ($G_2=\Spin9$
or $\Sp n\Sp1$ ($n\geq2$)). 

\subsubsection{} \label{3.1.2}
We address the case in which $\rho_1=\SO 3$. 
Formula \eqref{eq real tensor cohomo 1} gives
\begin{eqnarray*}
c(\rho)&=& c(\rho_2) + c(\SO 2,[c(\rho_2)-1]\RR^2) + c(H_2,2U_2)\\
&=&3c(\rho_2)-3+c(H_2,2U_2).
\end{eqnarray*}
We immediately get $c(\rho_2)\leq3$. If $c(\rho_2)=3$ then 
$c(H_2,2U_2)=2$, which is impossible by classification. 
If $c(\rho_2)=2$ then $c(H_2,2U_2)\leq 5$ 
implying $c(H_2,U_2)\leq 2$ which gives $c(H_2,2U_2)>5$, by 
classification. It follows that $c(\rho_2)=1$. 
We get two polar representations with $c=3$ ($\rho_2=\SO n$ or $\Spin7$),
one representation with $c=4$ ($\rho_2=\G)$,
a new example, namely, $\SO3\otimes\Sp n\Sp1$ ($n\geq2$) with $c=8$,
and $\SO3\otimes\Spin9$ with $c=9$. 

\subsubsection{}
For a representation $\rho_1\subset \SO n$, $n\geq 3$, 
we have $c(\rho_1\otimes_{\RR} \rho_2) \geq c(\SO 3\otimes \rho_2)$ and, thus, 
$\rho_2$ is of cohomogeneity one as follows from \ref{3.1.2}. 
In fact, we only need to deal with real tensor representations where both factors are of cohomogeneity one
and, more precisely, these are exceptional or one of $\SO n$, $\SU n$, $\U n$, or $\Sp1 \Sp2$. 
%We are only interested in irreducible representations.
Again, we rely on Lemma \ref{lem Tensors of cohom 1 representations} for the computation of the exact cohomogeneity in most cases. 
For $\rho_1=\SO n$, $n\geq 4$, we have 
$c(\SO n \otimes_{\RR} \U 2)= 6$ and $c(\SO n \otimes_{\RR} \SU 2)= 7$.
Among the exceptional cohomogeneity one representations as one of the factors we retrieve 
$c(\SO 4 \otimes_{\RR} \G)= 8$, $c(\SO 5 \otimes_{\RR} \G)\geq 11$, $c(\SO 4 \otimes_{\RR} \Spin 7)=5$ and  
$c(\SO 5 \otimes_{\RR} \Spin 7)\geq 9$. Moreover, we consider subrepresentations of the previous
when the cohomogeneity is smaller than $9$. We discard 
$c(\SU 2 \otimes_{\RR} \G) \geq c(\U 2 \otimes_{\RR} \G)\geq 10$, and get $c(\SU 2 \otimes_{\RR} \Spin 7)= 8$ 
and $c(\U 2 \otimes_{\RR} \Spin 7)= 7$.

\subsection{Quaternionic tensors}
Let us consider the tensor over the quaternionic algebra
of two simple Lie group quaternionic type representations.

\subsubsection{}
We begin assuming that $\rho_1$ is the standard $\Sp 1$-representation on $\HH$ giving the quaternionic structure.

Let $\rho_2$ be $(\Sp 1, \HH^n)$, $n\geq 1$.
The initial cases $n=1,2$ give small values of $c$. In turn, for $n\geq 4$, 
we have $c(\rho)\geq 4.n-6 \geq 10$. 
We are left with the case $n=3$ so that $\rho_2=(\Sp 1,\HH^3) \cong (\SU 2, S^5(\CC^2))$. 
This gives rise to $(\Sp 1\times \Sp 1, \HH \otimes_{\HH} \HH^3)$ which is of cohomogeneity $6$. 
%the natural description of the representation space 
%as degree $5$ %complex 
%homogeneous polynomials in two variables shows that it is of cohomogeneity $6$.
%%%%%NOT ORBIT EQUIVALENT TO SECOND FACTOR ONLY, THEN COHOMO =6 OR 7. IN THE FORMER CASE WE SHOULD HAVE A DIAGONAL CIRCLE INSIDETHE PIG.
%%%%%BUT THEN THIS CIRCLE MUST FIX AN 8 (REAL) DIMENSIONAL SUBSPACE WHILE *ANY* CIRCLE WILL FIX AT MOST 4 DIMENSIONS.

Assume now that the second factor is given by a bigger group, so that its dimension is at least $8$. 
%while the first factor is still the standard $\Sp 1$. 
In this case, we have 
\[
2.dim(G)\geq dim(G)+ 11 \geq dim(V), 
\]
and the representation is either polar or, else, belongs to a finite list, 
as stated in \cite[Proposition II]{K2}.
Among the non-polar cases we retrieve a cohomogeneity 6 representation given by 
\mbox{$(\Sp 1\otimes_{\HH}\Spin {11}, \HH^{16})$}.

\subsubsection{}
Consider the case where $\rho_1$ correspondst to $\Sp 1$ acting in a non-standard fashion, 
i.e., $\rho_1=(\Sp 1,\HH^m)$, $m\geq 2$. 
Then $\rho_1\otimes \rho_2 \subset \Sp m\otimes \Sp n$ with both $m,n\geq 2$. 
If $m\geq 3$ then 
\[
c(\Sp 1\times \Sp n,\HH^m\otimes \HH^n)\geq c(\Sp n,m.\HH^n)-3\geq 12. 
\]
Hence we may only admit  $m=2$ and $n\geq 2$. 
More generally, 
\[
c(\Sp 1\times G_2, \HH^2\otimes V_2)\geq c(G_2, 2.V_2)-3 \geq 2.c(G_2, V_2)-2.\] 
Our bound on the cohomogeneity then imposes $c(G_2,V_2)\leq 5$. 
It is easy to consider all such possible $\rho_2$ since they correspond to simple Lie group representations
of quaternionic type. We are thus led to a degree $4$ representation of $\Sp 1$ 
or else the standard cohomogeneity one action of $\Sp n$.
The former induces a representation whose cohomogeneity is at least $10$
while the latter gives a family $(\Sp 1\times \Sp n,\HH^2\otimes_{\HH}\HH^n)$ of cohomogeneity $3$ representations.
%and also exhibit copolarity one.

\subsubsection{}
We still need to consider the case where both factors are given by a bigger group than $\Sp 1$.
Using the Monotonicity Lemma it suffices to exclude cases of the form $\Sp 2\otimes \rho_2$.
We may bound the cohomogeneity from below as 
\[
c(\Sp 2\otimes G_2, \HH^2\otimes_{\HH} V_2) \geq 2.c(G_2, V_2) + 1 - 10 .
\]
Therefore, if $c(\rho)\leq 8$, we have $c(\rho_2)\leq 8$. 
Then either $\rho_2$ is $\Sp n$, of cohomogeneity one, or else there are $4$ such representations 
%	\SU 6, \HH^{10}
%	\Sp 3, \HH^7
%	\Spin {12}, \CC^{32}
%	\E 7,	\CC^{56}
which are discarded by a basic dimensional bound. 
The argument is analogous for the first factor, thus leading to $\rho=\Sp m \otimes_{\HH} \Sp n$ which gives a polar family.

\subsection{Complex tensors} 
In dealing with the remaining cases, we shall retrieve the real representations underlying complex tensor product representations. 
We now have a representation induced by restriction of $\U{n_1}\otimes_{\CC} \SU {n_2}\otimes_{\CC} \cdots \otimes_{\CC}  \SU {n_k}$
according to the number ``k'' of simple factors. 
In fact, we may have at most two simple factors since 
assuming otherwise leads to a subrepresentation of 
$\SU{m}\otimes_{\CC} \SU n\otimes_{\CC} \U p$ whose cohomogeneity is at least $16$,
%$c(\SU{m}\otimes_{\CC} \SU n\otimes_{\CC} \U p)\geq 16$, 
from a natural bound for the initial case \mbox{$(m,n,p)=(2,3,3)$}. 

\subsubsection{}
Let us consider the case of a simple Lie group representation $\rho_1$ 
%of complex or quaternionic type and cohomogeneity 6, 7, 8 or 9, 
extended by a circle factor, so that $6\leq c(\rho_1)\leq 9$. %cohomogeneity 6, 7, 8 or 9. 
We are looking for the non-polar examples and, thus, may assume that the simple factor itself is not polar 
(as observed in \cite{GL_OnOrbitSpaces}).
% one may see \cite[Table 12.3.1]{GL_OnOrbitSpaces}. SEE DISCUSSION SECTION ON SIMPLE GROUPS
We obtain two families of examples, first those of cohomogeneity $6$ given by: 
$(\U 6, \Lambda^3\CC^6)$, 
$(\U 1\Sp 3,\Lambda^3\CC^6\ominus\CC)$, 
$(\U 1\Spin {12},\CC^{32})$ and 
$(\U 1\E7,\CC^{56})$,
and those of cohomogeneity $8$ corresponding to 
$(\U 1\Spin{11},\CC^{32})$, and 
$(\U 1\Sp1,\CC^6\cong S^5(\CC^2))$.

These examples share structural similarities, lying in-between the representations given by restriction to the simple factor, $\rho_1$,
and the extension to $\Sp1\otimes_{\HH} \rho_1$, which corresponds to an s-representation in the cases at hand.
Furthermore, their cohomogeneities are known and differ by $3$, implying that the ``middle'' representations cannot be orbit equivalent to $\rho_1$. 
Thus, the cohomogeneity of the previous two families, of the form $\U1\otimes\rho_1$, equals $c(\rho_1)-1$. 
%In fact, the former ones come from cohomogeneity $4$ representations $(\Sp 1\otimes_{\HH} G,\RR^{4n})$ which are well known and polar. 
%lie in between cohomogeneity $7$ quaternionic type representations of the form $(G,\HH^n)$ and
%representation $(\Sp 1\otimes_{\HH} G,\RR^{4n})$ which are well known and polar. 
%Analogously, the latter two cases lie in between the simple factor cohomogeneity $9$ representations of the form
%$(G,\HH^n)$ and the cohomogeneity $6$ representations obtained as quaternionic tensors $(\Sp 1\otimes_{\HH}G,\RR^{4n})$.
%It should be pointed out that w

We have made use of the fact that 
the quaternionic tensor product of a given quaternionic type representation $\rho_1$ 
with the standard $(\Sp 1,\HH)$ restricted to a circle $\U 1\subset \Sp 1$ 
gives a representation that is equivariantly isomorphic to that 
of the complex tensor of $\rho_1$ with the standard $(\U 1,\CC)$.

\subsubsection{}\label{3.3.1}
Let $\rho$ be $\rho_1\otimes \rho_2 \subseteq \U m \otimes_{\CC} \SU n$, with $m\geq 2$. 
Assume first that $\rho_1=(\U 2, \CC^2)$ so that
$c(\rho_1\otimes \rho_2) \geq 2.c(\rho_2) -3 $.
The second factor is then given by a complex representation of a simple Lie group of cohomogeneity at most 5,
which is necessarily of complex type (after Remark \ref{rem tips for tensors by types}). 
%We recall that these are easily found among simple Lie group representations, 
%see Section \ref{section simple lie groups}.
The tensors with such representations can be discarded using that: 
$$c(\U 2\otimes \rho_2)\geq 2.dim_{\RR}(V_2) - 4 - dim(G_2),$$
unless they have cohomogeneity one. For the latter, given all the restrictions, 
we retrieve the polar $\U m \times \SU n$ complex-Grassmannian s-representations.
If the first factor was $\rho_1=(\SU2,\CC^2)$ the cohomogeneity would only increase.

\subsubsection{}
For $m\geq 3$ we have that $c(\rho_1\otimes \rho_2)\geq c(\U m\otimes \rho_2)\geq c(\U 2\otimes \rho_2)$ and,
%Considering $\rho_1=\U m$, $m\geq 3$, increases the cohomogeneity when compared to the previous case and, 
%therefore, implies that the second factor has $c=1$. 
therefore, the second factor is also of cohomogeneity one, by the previous case. 
In fact, both factors need to be of cohomogeneity one, thus leading to 
$\U m\otimes \Sp n$, $\SU m\otimes \Sp n$ and $\SU m\otimes \SU n$.
We compute their cohomogeneity with the aid of Lemma \ref{lem Tensors of cohom 1 representations}. 
%This leads to $\SU{m}\otimes_{\CC}\SU n$, $\SU{m}\otimes_{\CC}\U n$, $\SU{m}\otimes_{\CC}\Sp n$, $\U{m}\otimes_{\CC}\Sp n$
\qed

\medskip

\begin{table}[ht]
\[\begin{array}{|l|l|l|}%l|l|}
\hline
G		&		V		 &	Condition	\\
\hline
\hline
\SO3 				& \RR^9 			& - 	\\
\SU 5 				& S^2\CC^5 			& - 	\\
\SU {10} 			& \Lambda^2\CC^{10} 		& - 	\\
\U 1 \times \SO 3		& \CC^5				& -		\\
\U 1 \times \Sp 3		& \Lambda^2\CC^6\ominus\CC 	& -		\\
\U 1 \times\Sp 3		& \Lambda^3\CC^6 \ominus \CC^6	& -		\\
\U 1 \times \SU 6		& \Lambda^3\CC^6		& -		\\
\U 1 \times\Spin {12}		& \CC^{32}			& -		\\
\U 1 \times \E 7		& \CC^{56}			& -		\\
% tensors where both factors are of cohomogeneity 1
%real tensors
\SO 3\times \SU 2	&\RR^3\otimes_{\RR} \CC^2		&	-	\\
\SO 3\times \SU n	&\RR^3\otimes_{\RR} \CC^n		&	n\geq 4	\\
\SO 3\times \U n	&\RR^3\otimes_{\RR} \CC^n		&	n\geq 3	\\
\SO n\times \U 2	&\RR^n\otimes_{\RR} \CC^2		&	n\geq 4	\\
%quat tensors
\Sp 1\times \Sp 1 		&\HH \otimes_{\HH} \HH^3 	& -		\\
\Sp 1\times \Spin {11}		&\HH\otimes_{\HH} \HH^{16} 	& -		\\
%complex tensors
\SU 5\times \SU 5	&\CC^5\otimes_{\CC} \CC^5		&	-	\\
\U 3\times \Sp n	&\CC^3\otimes_{\CC} \HH^n		&	n\geq 3	\\		
\U n\times \Sp 2	&\CC^n\otimes_{\CC} \HH^2		&	n\geq 4	\\
\SU n\times \Sp 2	&\CC^n\otimes_{\CC} \HH^2		&	n=3,~n\geq 5	\\
\hline
  \end{array}\]
\caption{Non-polar irreducible representations of \mbox{cohomogeneity $6$}.}% of non-simple Lie groups.} 
\label{table cohom 6}
  \end{table}

\begin{table}[ht]
\[\begin{array}{|l|l|l|l|}%l|l|}
\hline
G		&		V			&	Condition	\\
\hline
\hline
 \Sp3 		& \Lambda^3\CC^6\ominus\CC^6			& - 		\\%& 7 	& q	& \F/\Sp 3\Sp 1	\\
 \SU 6 		& \Lambda^3\CC^6 				& - 		\\%& 7 		& q	& \E6/\SU 6\SU 2	\\
 \SU 6 		& S^2\CC^6 					& -		\\%5\leq n\leq 8 \\%& n+1 	& c 	& \Sp n/\U n	\\
 \SU {12} 	& \Lambda^2\CC^{12} 				& -		\\%{12}& \frac n2+1 & c & \SO{2n}/\U n\\
 \Spin{12} 	& \CC^{32} 					& - 		\\%5& 7			& q 	& \E7 /\Spin {12}\SU 2	\\%\mbox{\textsl{spin}}\\
 \E7 		& \CC^{56} 					& -		\\%& 7				& q	& \E8/\E7\SU 2\\
\U 1 \times \SU 3	&\mathfrak{su}(3)^{\CC} 		& -		\\
\U 1 \times \F		&\CC^{26}				& -		\\
\U 2 \times \Spin 7	&\CC^2\otimes_{\RR} \RR^8		& -		\\
% tensors where both factors are of cohomogeneity 1
%real tensors
\SO 3\times \SU 3 	&\RR^3\otimes_{\RR} \CC^3		& -		\\
\SO n\times \SU 2 	&\RR^n\otimes_{\RR} \CC^2		& n\geq 4	\\
%complex tensors
\SU 6\times \SU 6	&\CC^6\otimes_{\CC} \CC^6		&	-	\\
\SU 3\times \Sp n	&\CC^3\otimes_{\CC} \HH^n		& n\geq 3	\\
\SU 4\times \Sp 2	&\CC^4\otimes_{\CC} \HH^2		& -		\\
\hline
  \end{array}\]
\caption{Non-polar irreducible representations of \mbox{cohomogeneity $7$}.}% of non-simple Lie groups.} 
\label{table cohom 7}
  \end{table}

\begin{table}[ht]
\[\begin{array}{|l|l|}%l|l|}
\hline
G	\hspace{3mm}	&		V	\hspace{1cm}			 \\
\hline
\hline
\SO3 				& \RR^{11} 			\\%n=9,11 		& n-3 		& r	& -	\\
\SU {14} 			& \Lambda^2\CC^{14} 		\\%10 \leq n\leq 16, \mbox{$n$ even}	& \frac n2+1 & c & \SO{2n}/\U n\\
\SU 7 				& S^2\CC^7 			\\%5\leq n\leq 8 	& n+1 	& c 	& \Sp n/\U n	\\
\U 1\times \SU 2		& \HH^3 			\\
\U 1 \times \Spin {11}		& \HH^{16} 			\\
%
% tensors where both factors are of cohomogeneity 1
\SU 7\times \SU 7		&\CC^7\otimes_{\CC} \CC^7	\\
\SU 2\times \Spin 7		&\CC^2\otimes_{\RR} \RR^8	\\
\SO 4\times \G			& \RR^4 \otimes_{\RR} \RR^7 	\\
\SO 3\times \Sp1\Sp2		&\RR^3\otimes_{\RR}\HH^2	\\
\hline
\end{array}\]
\caption{Non-polar irreducible representations of \mbox{cohomogeneity $8$}.}% of non-simple Lie groups.} 
\label{table cohom 8}
  \end{table}

\end{document}